\documentclass{amsart}
\usepackage{amsfonts}

\usepackage{graphicx}
\usepackage{amscd}
\usepackage{amsmath}

\theoremstyle{plain}
\newtheorem{theorem}{Theorem}[section]

\newtheorem{lemma}[theorem]{Lemma}
\newtheorem{remark}[theorem]{Remark}
\newtheorem{question}[theorem]{Question}
\newtheorem{proposition}[theorem]{Proposition}
\theoremstyle{definition}

\numberwithin{equation}{section}

\input{tcilatex}

\begin{document}
\title[A Duality Proof of Tchakaloff's Theorem]{A Duality Proof of Tchakaloff's Theorem}
\author{Ra\'{u}l E. Curto}
\address{Department of Mathematics, The University of Iowa, Iowa City, Iowa 52242}
\email{curto@math.uiowa.edu}
\author{Lawrence A. Fialkow}
\address{Department of Computer Science, State University of New York, New Paltz, NY
12561}
\email{fialkowl@newpaltz.edu}
\thanks{Research partially supported by NSF grants. The second-named author was also
partially supported by the State University of New York at New Paltz
Research and Creative Projects Award Program.}

\begin{abstract}
Tchakaloff's Theorem establishes the existence of a quadrature rule of
prescribed degree relative to a positive, compactly supported measure that
is absolutely continuous with respect to Lebesgue measure on $\mathbb{R}^{d}$%
. Subsequent extensions were obtained by Mysovskikh and by Putinar. We
provide new proofs and partial extensions of these results, based on duality
techniques utilized by Stochel. \ We also obtain new uniqueness criteria in
the Truncated Complex Moment Problem.
\end{abstract}

\maketitle

\section{Introduction\label{sec1introduction}}

Tchakaloff's Theorem \cite[Th\'eor\`eme II]{TCH} establishes the existence
of a quadrature rule of prescribed degree relative to a positive, compactly
supported measure that is absolutely continuous with respect to Lebesgue
measure on $\mathbb{R}^{d}$. An extension to the case when the support is
unbounded was subsequently obtained by Mysovskikh \cite{MYS}, and in \cite
{PUT}, M. Putinar generalized these results to arbitrary positive Borel
measures. In the present note we provide new proofs, and partial extensions,
of Putinar's results, based on duality techniques utilized by J. Stochel 
\cite{STO} in a recent study of multivariable moment problems. We also
obtain some new uniqueness criteria in the Truncated Complex Moment Problem
(cf.\ \cite{TCMP1}, \cite{TCMP2}, \cite{TCMP3}).

For $t\equiv (t_{1},\dots ,t_{d})\in \mathbb{R}^{d}$, and for a multi-index $%
i\equiv (i_{1},\dots ,i_{d})\in \mathbb{Z}_{+}^{d}$, let $%
t^{i}=t_{1}^{i_{1}}\cdots t_{d}^{i_{d}}$ and $\left| i\right| =i_{1}+\cdots
+i_{d}$. Let $\mathbb{R}_{m,d}[t]\equiv \mathbb{R}_{m}[t_{1},\dots ,t_{d}]$
denote the space of real polynomials of total degree at most $m$ in $%
t_{1},\dots ,t_{d}$, and let $N_{m,d}:=\dim \mathbb{R}_{m,d}[t]$. As a
notational convenience, to indicate that each function in $\mathbb{R}%
_{m,d}[t]$ is (absolutely) integrable with respect to a positive Borel
measure $\mu $ on $\mathbb{R}^{d}$, we write $\mathbb{R}_{m,d}[t]\subseteq
L^{1}(\mu )$; note that in this case, the canonical map from $\mathbb{R}%
_{m,d}[t]$ into $L^{1}(\mu )$ need not be one-to-one. For a measure $\mu $
with closed support $K\subseteq \mathbb{R}^{d}$ and satisfying $\mathbb{R}%
_{m,d}[t]\subseteq L^{1}(\mu )$, a \textit{quadrature rule }of\textit{\
precision} (or \textit{degree}) $m$ and \textit{size} $N\;(<\infty )$
consists of \textit{nodes} $x_{1},\dots ,x_{N}$ in $K$ and positive weights $%
\rho _{1},\dots ,\rho _{N}$ such that 
\begin{equation}
\int p(t)\,d\mu (t)=\sum_{k=1}^{N}\rho _{k}p(x_{k})\qquad (p\in \mathbb{R}%
_{m,d}[t]).  \label{eq11}
\end{equation}

Putinar's first generalization of Tchakaloff's Theorem concerns the case of
compact support.

\begin{theorem}
\label{theorem1.1}\cite[Theorem 1]{PUT} Let $\mu $ be a positive, finite
Borel measure with compact support in $\mathbb{R}^{d}$. Given $m>0$, there
exists a quadrature rule for $\mu $ of degree $m$ with size $\leq N_{m,d}$.
\end{theorem}

\begin{remark}
For a positive Borel measure $\mu $ on $\mathbb{R}^{d}$ having convergent
moments up to at least degree $n$, let $N_{n,d;\mu }:=\dim \left\{ p|_{%
\limfunc{supp}\mu }:p\in \mathbb{R}_{n,d}[t]\right\} $. \ We show in Section
3 that the estimate\ $N_{m,d}$ in Theorem \ref{theorem1.1} can be refined to 
$N_{m,d;\mu }$ (cf. Theorem \ref{theorem3.5}). \ \ 
\end{remark}

For the case of unbounded support, Putinar's quadrature result assumes a
somewhat different form.

\begin{theorem}
\label{theorem1.2}\cite[Theorem 2]{PUT} Let $\mu $ be a positive Borel
measure supported in $\mathbb{R}^{d}$, and suppose that $\mu $ has
convergent moments up to at least degree $2m$, i.e., $\mathbb{R}%
_{2m,d}[t]\subseteq L^{1}(\mu )$. Then there exists a quadrature rule for $%
\mu $ of degree $2m-1$ with size $\leq N_{2m,d}$.
\end{theorem}

In \cite[Corollary 1]{PUT}, Putinar also obtained an analogue of Theorem \ref
{theorem1.2} valid in any degree (even or odd), provided $\limfunc{supp}\mu $
is contained in a proper closed convex cone of $\mathbb{R}^{d}$. In our main
result, which follows, we are able to remove the latter constraint on $%
\limfunc{supp}\mu $, so as to treat the ``even'' and ``odd'' cases together,
with an improved estimate for the size of a quadrature rule.

\begin{theorem}
\label{theorem1.3}(Generalized Tchakaloff Theorem, real case) Let $\mu $ be
a positive Borel measure on $\mathbb{R}^{d}$ having convergent moments up to
at least degree $n$. Then there exists a quadrature rule for $\mu $ of
degree $n-1$ with size $\leq 1+N_{n-1,d;\mu }$.
\end{theorem}

\begin{remark}
For a quadrature rule of degree $n-1$, the size estimate $1+N_{n-1,d;\mu }$
of Theorem \ref{theorem1.3} compares favorably with the estimate $N_{n,d}$
of Theorem \ref{theorem1.2} and \cite[Corollary 1]{PUT}. For the case of
compact support, Tchakaloff's size estimate $N_{n-1,d}$ (for a rule of
degree $n-1$) is known to be sharp \cite[pp. 131-133]{TCH}. It is thus
plausible that $1+N_{n-1,d;\mu }$ is sharp in the non-compact case, though
we do not have an example confirming this. Of course, for certain sets and
measures there exist Gaussian-type quadrature rules for which the size is
much smaller than that guaranteed by the above estimates \cite{MOL}, \cite
{RAL}, \cite{STR}, \cite{XU}. In particular, for $d=1$, any positive Borel
measure on $\mathbb{R}$, $[a,b]$, or $[0,+\infty )$, having moments up to
degree $n$, admits a Gaussian-type quadrature rule of degree $n$ with size $%
\leq \lbrack \frac{n}{2}]+1$ \cite{FIA3}.
\end{remark}

As we next describe, Theorems \ref{theorem1.1} - \ref{theorem1.3} are
closely related to the following Truncated Multivariable Moment Problem
(TMMP) \cite[Chapter 7]{TCMP1}. For a real multisequence $\beta \equiv \beta
^{(m)}=\{\beta _{i}\}_{i\in \mathbb{Z}_{+}^{d},\text{ }\left| i\right| \leq
m}$, TMMP seeks to characterize the existence of a positive Borel measure $%
\mu $ supported in $\mathbb{R}^{d}$ such that 
\begin{equation}
\beta _{i}=\int t^{i}\,d\mu (t)\qquad (\left| i\right| \leq m);  \label{eq12}
\end{equation}
a measure $\mu $ as in (\ref{eq12}) is said to be a \textit{representing
measure} for $\beta $. The $K$\textit{-moment problem} further requires $%
\limfunc{supp}\mu \subseteq K$, where $K$ is a prescribed closed subset of $%
\mathbb{R}^{d}$. The following basic question remains open.

\begin{question}
\label{question1.3}If $\beta \equiv \beta ^{(m)}$ has a representing
measure, does $\beta $ admit a \emph{finitely atomic} representing measure $%
( $i.e., a measure of the form $\sum_{k=1}^{N}\rho _{k}\delta _{x_{k}}$,
where $1\leq N<\infty $, each $\rho _{k}>0$, and $\delta _{x_{k}}$ is the
point mass at $x_{k}\in \mathbb{R}^{d})$? More generally, if $\beta $ has a
representing measure supported in $K$, does $\beta $ have a finitely atomic
representing measure supported in $K$?
\end{question}

Let $\mu $ be a representing measure for $\beta \equiv \beta ^{(m)}$ and let 
$K:=\limfunc{supp}\mu $. The existence of a finitely atomic representing
measure in the $K$-moment problem for $\beta $ is equivalent to the
existence of a quadrature rule for $\mu $ of degree $m$. For $d=1$ and $K=%
\mathbb{R}$, $[a,b]$, or $[0,\infty )$, \cite{FIA3} implies that Question 
\ref{question1.3} has an affirmative answer relative to representing
measures in the $K$-moment problem. For $d\geq 1$, Theorem \ref{theorem1.1}
shows that if $\beta $ has a representing measure supported in a compact set 
$K$, then $\beta $ admits a finitely atomic representing measure supported
in $K$. Theorem \ref{theorem1.3} implies that if $\beta \equiv \beta ^{(m)}$
has a representing measure $\mu $ satisfying $\mathbb{R}_{m+1,d}[t]\subseteq
L^{1}(\mu )$ and $\limfunc{supp}\mu \subseteq K$ (closed), then $\beta $
admits a finitely atomic representing measure supported in $K$.

The proof of Theorem \ref{theorem1.2} in \cite{PUT} rests on convexity
arguments similar to those in \cite{TCH}, and also on some rather subtle
convergence arguments. Our proof of Theorem \ref{theorem1.3} is based on a
different approach. Given a multisequence $\beta \equiv \beta ^{(m)}$ and a
prescribed closed set $K\subseteq \mathbb{R}^{d}$, let $M_{K}(\beta )$
denote the set of representing measures for $\beta $ supported in $K$. In
Proposition \ref{prop2.3}, we show that any extreme point of $M_{K}(\beta )$
is finitely atomic. This result is based on a technique of J.P. Gabardo \cite
{GAB}, who attributes the underlying idea to Naimark (cf. \cite{AKH}).
Proposition \ref{prop2.3} yields some insight into Question \ref{question1.3}%
; in particular, if there is a unique representing measure for $\beta $, it
is finitely atomic. Let $\mu $ be a measure satisfying the hypothesis of
Theorem \ref{theorem1.3}, i.e., $\mathbb{R}_{n,d}[t]\subseteq L^{1}(\mu )$.
In Section \ref{sec3dualityproof} we employ duality results of Stochel \cite
{STO} to show that a certain convex set of representing measures for $\beta
^{(n-1)}[\mu ]$ admits an extreme point, which then acts as a quadrature
rule for $\mu $ of degree $n-1$.

In Section \ref{sec4unique} we present some uniqueness results in the
following Truncated Complex Moment Problem (TCMP) \cite{TCMP1}. Let $\gamma
\equiv \gamma ^{(2n)}:\gamma _{00},$ $\gamma _{01},$ $\gamma _{10},$ $\gamma
_{02},$ $\gamma _{11},$ $\gamma _{20},\dots ,$ $\gamma _{0,2n},\dots ,$ $%
\gamma _{2n,0}$ denote a sequence of complex numbers. TCMP entails
characterizing the existence of a positive Borel measure $\mu $ on $\mathbb{C%
}$ such that 
\begin{equation*}
\gamma _{ij}=\int \bar{z}^{i}z^{j}\,d\mu (z)\qquad (0\leq i+j\leq 2n).
\end{equation*}
We encode the data $\gamma $ in a moment matrix $M(n)(\gamma )$ 
\cite[Chapter 2]{TCMP1}, whose successive columns are labeled $1,Z,\bar{Z}%
,\dots ,Z^{n},\dots ,\bar{Z}^{n}$. If $\gamma $ has a representing measure,
then $M(n)(\gamma )$ is positive semidefinite and \textit{recursively
generated} (cf. \cite[(3.2) and Remark 3.15(ii)]{TCMP1}); moreover, if $\mu $
is a representing measure for $\gamma $, then $\limfunc{card}\limfunc{supp}%
\mu \geq \limfunc{rank}M(n)$ \cite[Corollary 3.7]{TCMP1}. In \cite[Chapter 5]
{TCMP1}, for the case of \textit{flat data}, where $M(n)\geq 0$ and $%
\limfunc{rank}$ $M(n)=\limfunc{rank}M(n-1)$, we established the existence of
(and explicitly constructed) a unique finitely atomic representing measure.
In Proposition \ref{prop4.1} we show that this measure (with $\limfunc{rank}%
M(n)$ atoms) is actually the unique representing measure for $\gamma $.

We say that $M(n)(\gamma )$ admits an \textit{analytic relation} if there
exist $k\leq n$ and scalars $a_{ij}\in \mathbb{C}$ $(0\leq i+j<k)$ such that
in the column space of $M(n)$ there is a dependence relation 
\begin{equation}
Z^{k}=\sum_{0\leq i+j<k}a_{ij}\bar{Z}^{i}Z^{j}.  \label{equ13}
\end{equation}
In \cite[Theorem 3.1]{TCMP2} we proved that if $M(n)$ is positive and
recursively generated, and if $M(n)$ admits an analytic relation with $k\leq %
\left[ \frac{n}{2}\right] +1$, then $\gamma $ admits a unique finitely
atomic representing measure (with $\limfunc{rank}M(n)$ atoms); by contrast,
for $k>\left[ \frac{n}{2}\right] +1$, there need not be any representing
measure \cite[Example 2.14]{TCMP3}, or there may be a finitely atomic
representing measure, but none with as few as $\limfunc{rank}M(n)$ atoms 
\cite[Theorem 3.1]{FIA1}. In Proposition \ref{prop4.2} we prove that if $%
\gamma $ has a representing measure and $M(n)(\gamma )$ admits an analytic
relation (as in (\ref{equ13})) for some $k\leq n$, then $\gamma $ has a
unique representing measure, which is finitely atomic, with at most $k^{2}$
atoms. \ This result depends on the fact that a polynomial of the form $%
z^{k}-q(z,\bar{z})\;(\deg q<k)$ has at most $k^{2}$ roots (Proposition \ref
{lemma4.3}).

\section{Extreme points and finitely atomic representing measures\label%
{sec2extremepoints}}

For $z\equiv (z_{1},\dots ,z_{d})\in \mathbb{C}^{d}$ and $j\equiv
(j_{1},\dots ,j_{d})\in \mathbb{Z}_{+}^{d}$, let $z^{j}=z_{1}^{j_{1}}\cdots
z_{d}^{j_{d}}$ and let $\left| j\right| =j_{1}+\cdots +j_{d}$. Let $\mathbb{C%
}_{m,d}[z,\bar{z}]$ denote the complex polynomials $p(z,\bar{z})$ of total
degree at most $m$. For a closed set $K\subseteq \mathbb{C}^{d}$, $\mathbb{C}%
_{m,d}[z,\bar{z}]\mid _{K}$ denotes the vector space of restrictions to $K$
of polynomials in $\mathbb{C}_{m,d}[z,\bar{z}]$. Given a complex sequence $%
\gamma \equiv \gamma ^{(m)}=\{\gamma _{ij}\}_{i,j\in \mathbb{Z}_{+}^{d},%
\text{ }\left| i\right| +\left| j\right| \leq m}$, and a closed set $%
K\subseteq \mathbb{C}^{d}$, the Multivariable Truncated Complex $K$-Moment
Problem entails characterizing the existence of a positive Borel measure $%
\mu $, supported in $K$, such that 
\begin{equation*}
\gamma _{ij}=\int \bar{z}^{i}z^{j}\,d\mu ,\qquad \left| i\right| +\left|
j\right| \leq m
\end{equation*}
(where $\bar{z}^{i}=\bar{z}_{1}^{i_{1}}\cdots \bar{z}_{d}^{i_{d}}$). The
Multivariable Full Complex $K$-Moment Problem in $\mathbb{C}^{d}$ concerns
the analogous problem, for a sequence $\gamma ^{(\infty )}\equiv \{\gamma
_{i,j}\}_{i,j\in \mathbb{Z}_{+}^{d}}$, which prescribes moments of all
orders. The Full and Truncated Complex $K$-Moment Problems on $\mathbb{C}%
^{d} $ are equivalent, respectively, to corresponding moment problems on $%
\mathbb{R}^{2d}$ (cf.\ \cite[Chapters 6 and 7]{TCMP1}, \cite[Section 5]
{TCMP4}). The following result of J. Stochel provides the connection between
the Full and Truncated Multivariable Complex $K$-Moment Problems.

\begin{theorem}
\label{theorem2.1}\cite[Theorem 4]{STO} $\gamma ^{(\infty )}$ has a
representing measure supported in a closed set $K\subseteq \mathbb{C}^{d}$
if and only if, for each $m>0$, $\gamma ^{(m)}$ admits a representing
measure supported in $K$.
\end{theorem}

Assume that $\gamma \equiv \gamma ^{(m)}$ admits a representing measure
supported in a closed set $K\subseteq \mathbb{C}^{d}$. Consider the convex
set 
\begin{equation*}
M_{K}(\gamma )=\{\nu :\nu \text{ is a representing measure for }\gamma \text{
and }\limfunc{supp}\nu \subseteq K\}.
\end{equation*}
It is not known whether $M_{K}(\gamma )$ always has an extreme point. The
proof of the following result is motivated by an argument in 
\cite[Proposition 2.5 and Corollary 2.6]{GAB}. For $\nu \in M_{K}(\gamma )$,
let $\mathbb{C}_{m,d}[z,\bar{z}](\nu )$ denote the image of $\mathbb{C}%
_{m,d}[z,\bar{z}]$ in $L^{1}(\nu )$ under the canonical projection.

\begin{proposition}
\label{prop2.2}If $M_{K}(\gamma )$ has an extreme point $\nu $, then $\nu $
is finitely atomic, with $\limfunc{card}\limfunc{supp}\nu \leq \dim \mathbb{C%
}_{m,d}[z,\bar{z}](\nu )$ ($\leq \dim \mathbb{C}_{m,d}[z,\bar{z}]\;|_{K}$).
\end{proposition}

\begin{proof}
Since $\nu $ is a representing measure for $\gamma ^{(m)}$, we may consider $%
\mathcal{L}:=\mathbb{C}_{m,d}[z,\bar{z}](\nu )\subseteq L^{1}(\nu )$, and we
claim that $\mathcal{L}$ is dense in $L^{1}(\nu )$.

Since $\nu $ is finite, it follows that $L^{1}(\nu )^{\ast }=L^{\infty }(\nu
)$ \cite[Theorem III.5.6]{CON}; thus, if $\mathcal{L}$ is not dense, there
exists $f\in L^{\infty }(\nu )$, $f\neq 0$, such that 
\begin{equation}
\int pf\,d\nu =0\qquad (p\in \mathcal{L}).  \label{eq2.1}
\end{equation}
Since $\nu \geq 0$ and $\mathcal{L}$ is self-adjoint, we may replace $f$ by $%
\frac{1}{2\left\| f+\bar{f}\right\| _{\infty }}$\thinspace $(f+\bar{f})$,
and we may thus assume that $f$ is real, with $\left\| f\right\| _{\infty
}\leq \frac{1}{2}$. Thus $\nu _{1}:=(1+f)\nu $ and $\nu _{2}:=(1-f)\nu $ are
positive, and (\ref{eq2.1}) implies that they belong to $M_{K}(\gamma )$.
Since $\nu =\frac{1}{2}\nu _{1}+\frac{1}{2}\nu _{2}$, we have a
contradiction to the hypothesis that $\nu $ is an extreme point of $%
M_{K}(\gamma )$.

Now $\mathcal{L}$ is dense in $L^{1}(\nu )$, and since $\mathcal{L}$ is
finite dimensional, we have $\mathcal{L}=L^{1}(\nu )$, whence $r:=\dim 
\mathcal{L}=\dim L^{1}(\nu )=\dim L^{1}(\nu )^{\ast }=\dim L^{\infty }(\nu )$%
. Suppose $\limfunc{supp}\nu $ contains distinct points $z_{1},\dots
,z_{r+1} $, and let $\{D_{i}\}_{i=1}^{r+1}$ denote mutually disjoint closed
disks of positive radii such that $z_{i}\in D_{i}$ $(1\leq i\leq r+1)$. Then 
$\{\chi _{D_{i}}\}_{i=1}^{r+1}$ is linearly independent in $L^{\infty }(\nu
) $; this contradiction implies 
\begin{eqnarray*}
\limfunc{card}\limfunc{supp}\nu &\leq &r=\dim \mathbb{C}_{m,d}[z,\bar{z}%
](\nu ) \\
&\leq &\dim \mathbb{C}_{m,d}[z,\bar{z}]\mid _{\limfunc{supp}\nu }\leq \dim 
\mathbb{C}_{m,d}[z,\bar{z}]\mid _{K}.
\end{eqnarray*}
\end{proof}

The analogue of Proposition \ref{prop2.2} for a real multisequence $%
\beta^{(m)}$ and a closed subset $K\subseteq\mathbb{R}^{d}$ can be proved by
a straightforward modification of the preceding argument; we omit the
details.

\begin{proposition}
\label{prop2.3}If $M_{K}(\beta ^{(m)})$ has an extreme point $\nu $, then $%
\nu $ is finitely atomic, with $\limfunc{card}\limfunc{supp}\nu \leq \dim 
\mathbb{R}_{m,d}[t](\nu )\;(\leq \dim \mathbb{R}_{m,d}[t]\mid _{K})$.
\end{proposition}

\section{A duality proof of Tchakaloff's Theorem\label{sec3dualityproof}}

The main result of this section is the following complex version of the
generalized Tchakaloff theorem; we treat the complex case first mostly as a
convenience, since the tools we require from \cite{STO} are formulated in
terms of the complex moment problem.

For a positive Borel measure $\mu $ on $\mathbb{C}^{d}$ having convergent
moments up to at least degree $n$, recall that $\mathbb{C}_{n,d}[z,\bar{z}%
](\mu )$ denotes the image of $\mathbb{C}_{n,d}[z,\bar{z}]$ in $L^{1}(\mu )$
under the canonical projection, and let $\mathcal{N}_{n,d}(\mu ):=\dim 
\mathbb{C}_{n,d}[z,\bar{z}](\mu )$ and $N_{n,d;\mu }:=\dim \mathbb{C}%
_{n,d}[z,\bar{z}]\;|_{\limfunc{supp}\mu }$.

\begin{theorem}
\label{theorem3.1}(Generalized Tchakaloff Theorem, complex case). Let $\mu $
be a positive Borel measure on $\mathbb{C}^{d}$ having convergent moments up
to at least degree $n$, and let $K:=\limfunc{supp}\mu $. Then, for some $%
N\leq 1+N_{n-1,d;\mu }$, there exist nodes $z_{1},\dots ,z_{N}\in K$, and
positive weights $\rho _{1},\dots ,\rho _{N}$, such that 
\begin{equation*}
\int_{K}p(z,\bar{z})\,d\mu (z)=\sum_{k=1}^{N}\rho _{k}p(z_{k},\bar{z}%
_{k})\qquad (p\in \mathbb{C}_{n-1,d}[z,\bar{z}]).
\end{equation*}
\end{theorem}

To prove Theorem \ref{theorem3.1}, we first introduce some preliminary
results and notation concerning duality. Let $X$ be a locally compact
Hausdorff space. A continuous function $f:X\rightarrow \mathbb{C}$ \textit{%
vanishes at infinity} if, for each $\epsilon $ $>0$, there is a compact set $%
C_{\epsilon }\subseteq X$ such that $X\setminus C_{\epsilon }\subseteq
\{x\in X:\left| f(x)\right| <$ $\epsilon \}$. Let $C_{0}(X)$ denote the
Banach space of all functions on $X$ which vanish at infinity, equipped with
the norm $\left\| f\right\| _{\infty }:={\sup_{x\in X}}\left| f(x)\right| $.
The space $C_{c}(X)$ of continuous functions with compact support, is norm
dense in $C_{0}(X)$ \cite[III.1, Exercise 13]{CON}; when $X$ is compact, $%
C_{0}(X)=C_{c}(X)=C(X)$, the space of continuous complex-valued functions on 
$X$. The Riesz Representation Theorem \cite[C.18]{CON} states that $%
C_{0}(X)^{\ast }$, the dual space of $C_{0}(X)$, is isometrically isomorphic
to $M(X)$, the space of finite regular complex Borel measures on $X$
(equipped with the norm $\left\| \mu \right\| :=\left| \mu (X)\right| $);
under this duality, corresponding to $\mu \in M(X)$ is the functional $\hat{%
\mu}$ on $C_{0}(X)$ defined by $\hat{\mu}(f):=\int f\,d\mu $.

We now focus on the case where $\mu $ is a positive Borel measure on $%
\mathbb{C}^{d}$ with convergent moments up to (at least) order $n$, i.e., $%
\mathbb{C}_{n,d}[z,\bar{z}]\subseteq L^{1}(\mu )$. Let $K:=\limfunc{supp}\mu 
$; without loss of generality, in the sequel we normalize $\mu $ so that $%
\mu (K)=1$. Since the monomials $\bar{z}^{i}z^{j}\;(\left| i\right| +\left|
j\right| \leq n)$ are absolutely integrable, we may consider the
corresponding moments of $\mu $, 
\begin{equation*}
\gamma _{ij}:=\int_{K}\bar{z}^{i}z^{j}\,d\mu (z),\qquad \left| i\right|
+\left| j\right| \leq n.
\end{equation*}
Let $\rho _{0}(z):=\left\| z\right\| ^{n}$ (where, as usual, $\left\|
z\right\| :=(\left| z_{1}\right| ^{2}+\cdots +\left| z_{d}\right|
^{2})^{1/2} $). Let $\Gamma :=\int_{K}\rho _{0}(z)\,d\mu (z)$. If $n$ is
even, say $n=2m$, then $\left\| z\right\| ^{n}=(\bar{z}_{1}z_{1}+...+\bar{z}%
_{d}z_{d})^{m}\in \mathbb{C}_{n,d}[z,\bar{z}]$, so $\Gamma <\infty $. For
the case when $n$ is odd, to see that $\Gamma <\infty $ note that by the
equivalence of all norms on $\mathbb{C}^{d}$, there exists a constant $M>0$
such that for every $z\in \mathbb{C}^{d}$, $\left\| z\right\| \leq M(\left|
z_{1}\right| +...+\left| z_{d}\right| )$. Then $\int_{K}\left\| z\right\|
^{n}\,d\mu (z)\leq M^{n}\int_{K}(\left| z_{1}\right| +...+\left|
z_{d}\right| )^{n}\,d\mu (z)$, and the latter integral is convergent, since
for every multi-index $i$ with $\left| i\right| =n$, $z^{i}$ ($%
=z_{1}^{i_{1}}\cdot ...\cdot z_{d}^{i_{d}}$) is absolutely integrable.

Let $\gamma \equiv \gamma ^{(n)}[\mu ]:=\{\gamma _{i,j}\}_{0\leq \left|
i\right| +\left| j\right| \leq n}$, and set $V(\mu ;n):=\left\{ \nu \in
M(K):\nu \geq 0\right. ,\;\gamma _{ij}=\int_{K}\bar{z}^{i}z^{j}\,d\nu
(z),\left| i\right| +\left| j\right| \leq n-1,\;\int_{K}\left\| z\right\|
^{n}\,d\nu (z)\leq \Gamma \}.$ Observe that $V(\mu ;n)$ is convex, and is
nonempty since $\mu \in V(\mu ;n)$. For $\nu \in V(\mu ;n)$, $\left\| \nu
\right\| =\nu (K)=\gamma _{00}=1$, so $V(\mu ;n)$ embeds as a subset of $%
B_{1}(C_{0}(K)^{\ast })$, the closed unit ball of $C_{0}(K)^{\ast }$; recall
that $B_{1}(C_{0}(K)^{\ast })$ is weak-$\ast $ compact and metrizable 
\cite[Theorems V.3.1 and V.5.1]{CON}.

\begin{proposition}
\label{prop3.2}$V(\mu;n)$ is weak-$*$ closed in $B_{1}(C_{0}(K)^{*})$.
\end{proposition}

To prove Proposition \ref{prop3.2} we rely on the following technical result
of Stochel.

\begin{proposition}
\label{prop3.3}(\cite[Proposition 1]{STO}) Let $F$ be a nonempty closed
subset of $\mathbb{C}^{d}$ and let $\rho $ be a non-negative continuous
function on $F$. Assume that $\{\mu _{w}\}_{w\in \Omega }$ is a net of
finite positive Borel measures on $F$ and $\mu $ is a finite positive Borel
measure on $F$ such that

\begin{enumerate}
\item[$\mathrm{(i)}$]  ${\lim_{w\in \Omega }}\int_{F}f\,d\mu
_{w}=\int_{F}f\,d\mu $ $(f\in C_{c}(F))$, and$\rule[-15pt]{0pt}{15pt}$

\item[$\mathrm{(ii)}$]  ${\sup_{w\in \Omega }}\int_{F}\rho \,d\mu
_{w}<+\,\infty $.
\end{enumerate}

Then $\int_{F}\rho \,d\mu \leq {\sup_{w\in \Omega }}\int_{F}\rho \,d\mu _{w}$%
, and $\int_{F}f\rho \,d\mu ={\lim_{w\in \Omega }}\int f\rho \,d\mu _{w}$ $%
(f\in C_{0}(F))$. Moreover, if the set $\{z\in F:\rho (z)\leq r\}$ is
compact for some $r>0$, then $\int_{F}f\,d\mu ={\lim_{w\in \Omega }}%
\int_{F}f\,d\mu _{w}$ for every $f:F\rightarrow \mathbb{C}$ such that $\frac{%
f}{1+\rho }\in C_{0}(F)$.
\end{proposition}

\noindent(The estimate on $\int_{F}\rho\,d\mu$ is not part of the statement
of \cite[Proposition 1]{STO}, but is established in the proof.)

\begin{proof}[Proof of Proposition $\mathrm{3.2}$]
Since $B_{1}(C_{0}(K)^{\ast})$ is weak-$\ast$ metrizable, to establish that $%
V(\mu;n)$ is weak-$\ast$ closed, it suffices to show that it is closed under
limits of sequences. Let $\{\nu _{k}\}_{k=1}^{\infty}\subseteq V(\mu;n)$ and
suppose $\{\nu_{k}\}$ is weak-$\ast$ convergent to $\Lambda\in
B_{1}(C_{0}(K)^{\ast})$, i.e., $\Lambda(f)=\lim_{k\rightarrow\infty}\int
f\,d\nu_{k}$ $(f\in C_{0}(K))$. Clearly, $\Lambda\geq0$, since each $%
\nu_{k}\geq0$. Thus, by the Riesz Representation Theorem, there exists a
finite, positive Borel measure $\nu$, $\limfunc{supp}\nu\subseteq K$, such
that $\Lambda=\hat{\nu}$.

We claim that $\nu \in V(\mu ;n)$. Let $\rho =\rho _{0}$, i.e., $\rho
(z):=\left\| z\right\| ^{n}(z\in K)$. Since $\nu _{k}\in V(\mu ;n)$, then $%
\int_{K}\rho (z)\,d\nu _{k}(z)\leq \Gamma $ $(\equiv \int_{K}\rho (z)\,d\mu
) $, whence Proposition \ref{prop3.3} implies $\int_{K}\rho (z)\,d\nu \leq
\Gamma $. To complete the proof, we will show that for $\left| i\right|
+\left| j\right| \leq n-1$, $f_{ij}(z,\bar{z}):=\bar{z}^{i}z^{j}$ satisfies 
\begin{equation}
\int_{K}f_{ij}\,d\nu (z)=\gamma _{ij}.  \label{eq3.1}
\end{equation}
For each $k$, $\int_{K}f_{ij}\,d\nu _{k}(z)=\gamma _{ij}$; Proposition \ref
{prop3.3} implies that to establish (\ref{eq3.1}) it suffices to verify that 
$\frac{f_{ij}}{1+\rho }\in C_{0}(K)$ $(\left| i\right| +\left| j\right| \leq
n-1)$. Let $L>1$ and suppose $\left\| z\right\| ^{2}>L^{2}$, i.e., $z$ is in
the complement of the compact set $\{z\in \mathbb{C}^{d}:\left\| z\right\|
\leq L\}$. Choose $i(z),$ $1\leq i(z)\leq d$, such that $\left|
z_{i(z)}\right| \geq \left| z_{i}\right| ,$ $1\leq i\leq d$ ($i(z)$ depends
on $z$). Then 
\begin{equation*}
\left| z_{i(z)}\right| ^{2}\geq \frac{\left\| z\right\| ^{2}}{d}>\frac{L^{2}%
}{d},
\end{equation*}
so 
\begin{equation*}
\frac{1}{\left| z_{i(z)}\right| }<\frac{\sqrt{d}}{L}.
\end{equation*}
Now 
\begin{align*}
\frac{\left| f_{ij}\right| }{1+\rho }& =\frac{\left| z_{1}\right|
^{i_{1}+j_{1}}\cdots \left| z_{d}\right| ^{i_{d}+j_{d}}}{1+\left( \left|
z_{1}\right| ^{2}+\cdots +\left| z_{d}\right| ^{2}\right) ^{n/2}}\leq \frac{%
\left| z_{i(z)}\right| ^{\left| i\right| +\left| j\right| }}{\left|
z_{i(z)}\right| ^{n}} \\
& =\frac{1}{\left| z_{i(z)}\right| ^{n-\left| i\right| -\left| j\right| }}<%
\frac{d^{n}}{L}\longrightarrow 0\qquad (L\longrightarrow +\,\infty ).
\end{align*}
Thus $\frac{f_{ij}}{1+\rho }\in C_{0}(K)$ and the proof is complete.
\end{proof}

We next require a variant of Proposition \ref{prop2.2}.

\begin{lemma}
\label{lemma3.4}If $\nu $ is an extreme point of $V(\mu ;n)$, then $\nu $ is
finitely atomic, with $\limfunc{card}\limfunc{supp}\nu \leq 1+\dim \mathbb{C}%
_{n-1,d}[z,\bar{z}](\nu )$.
\end{lemma}

\begin{proof}
Since $\nu $ is a representing measure for $\{\gamma _{ij}\}_{\left|
i\right| +\left| j\right| \leq n-1}$, and $\int \left\| z\right\| ^{n}\,d\nu
<+\,\infty $, we may consider the subspace of $L^{1}(\nu )$ defined by $%
\mathcal{M}:=\{[p]+\alpha \lbrack \rho _{0}]:p\in \mathbb{C}_{n-1,d}[z,\bar{z%
}]$, $\alpha \in \mathbb{C}\}$. We claim that $\mathcal{M}$ is dense in $%
L^{1}(\nu )$. Since $\rho _{0}\geq 0$, $\mathcal{M}$ is self-adjoint; thus,
if $\mathcal{M}$ is not dense, it follows as in the proof of Proposition \ref
{prop2.2} that there exists $f:\limfunc{supp}\nu \rightarrow \mathbb{R}$, $%
\left\| f\right\| _{\infty }\leq \frac{1}{2}$, such that $\int pf\,d\nu =0$
for every $p\in \mathcal{M}$. Let $\nu _{1}:=(1+f)\nu $ and $\nu
_{2}:=(1-f)\nu $; for $i=1,2,\nu _{i}\geq 0$, $\nu _{i}$ is a representing
measure for $\{\gamma _{ij}\}_{\left| i\right| +\left| j\right| \leq n-1}$,
and 
\begin{equation*}
\int \left\| z\right\| ^{n}\,d\nu _{i}=\int \left\| z\right\| ^{n}\,d\nu \pm
\int \left\| z\right\| ^{n}f\,d\nu =\int \left\| z\right\| ^{n}\,d\nu \leq
\Gamma .
\end{equation*}
Thus $\nu _{i}\in V(\mu ;n)$ $(i=1,2)$, and since $\nu =\frac{1}{2}\nu _{1}+%
\frac{1}{2}\nu _{2}$, we have a contradiction to the hypothesis that $\nu $
is an extreme point. The rest of the proof is identical to that of
Proposition \ref{prop2.2}; in particular, since $\mathcal{M}$ is finite
dimensional, 
\begin{equation*}
\limfunc{card}\limfunc{supp}\nu \leq \dim L^{\infty }(\nu )=\dim L^{1}(\nu
)=\dim \mathcal{M}\leq 1+\dim \mathbb{C}_{n-1,d}[z,\bar{z}](\nu ).
\end{equation*}
\end{proof}

\begin{proof}[Proof of Theorem $\mathrm{3.1}$]
Since $\mu \in V(\mu ;n)$, Proposition \ref{prop3.2} implies that $V(\mu ;n)$
is a nonempty weak-$\ast $ compact convex subset of $B_{1}(C_{0}(K)^{\ast })$
(where $K:=\limfunc{supp}\mu $). It follows from the Krein-Millman Theorem
that $V(\mu ;n)$ has an extreme point, and Lemma \ref{lemma3.4} implies that
any such extreme point $\nu $ corresponds to a quadrature rule for $\mu $ of
degree $n-1$ with size at most $1+N_{n-1,d}(\nu )\leq 1+N_{n-1,d;\nu }\leq
1+N_{n-1,d;\mu }$.
\end{proof}

We next turn to the complex version of Theorem \ref{theorem1.1}.

\begin{theorem}
\label{theorem3.5}Suppose $\mu $ is a positive finite Borel measure on $%
\mathbb{C}^{d}$ with compact support $K$. Given $m>0$, there exist $N\leq
\dim (\mathbb{C}_{m,d}[z,\bar{z}]\mid _{K})$, nodes $z_{1},\dots ,z_{N}$ in $%
K$, and positive weights $\rho _{1},\dots ,\rho _{N}$ such that 
\begin{equation*}
\int_{K}p(z,\bar{z})\,d\mu (z)=\sum_{i=1}^{N}\rho _{i}p(z_{i},\bar{z}%
_{i})\qquad (p\in \mathbb{C}_{m,d}[z,\bar{z}]).
\end{equation*}
\end{theorem}

\begin{proof}
Since $C_{0}(K)=C(K)$, it is straightforward to modify the proof of
Proposition \ref{prop3.2} to conclude that $M_{K}(\gamma ^{(m)}[\mu ])$ is
weak-$\ast $ closed in $B_{1}(C_{0}(K)^{\ast })$. Thus $M_{K}(\gamma
^{(m)}[\mu ])$ has an extreme point $\nu $, and Proposition \ref{prop2.2}
implies $\limfunc{card}\limfunc{supp}\nu \leq \dim \mathbb{C}_{m,d}[z,\bar{z}%
]\mid _{\limfunc{supp}\nu }\leq \dim \mathbb{C}_{m,d}[z,\bar{z}]\mid _{K}$.
\end{proof}

\begin{proof}[Proofs of Theorems $\mathrm{1.1-}$ $\mathrm{1.3}$]
The proofs of Theorems \ref{theorem3.1} and \ref{theorem3.5} depend on
duality results for $L^{1}(\mu )^{\ast }$ and $C_{0}(K)^{\ast }$, where $\mu 
$ is a positive Borel measure on $\mathbb{C}^{d}$ and $K$ is a closed subset
of $\mathbb{C}^{d}$. These duality results, including those of Stochel
(e.g., Proposition \ref{prop3.3}), admit exact analogues for the real case,
where $\mathbb{C}^{d}$ is replaced by $\mathbb{R}^{d}$. Thus we may formally
repeat the proof of Theorem \ref{theorem3.1} to obtain Theorem \ref
{theorem1.3} (and Theorem \ref{theorem1.2}), and similarly for Theorem \ref
{theorem3.5} and Theorem \ref{theorem1.1}. (For the case when $d=2p$,
alternate proofs can be based on the equivalence of moment problems on $%
\mathbb{R}^{2p}$ with those on $\mathbb{C}^{p}$.)
\end{proof}

Implicit in \cite[pp. 127-129]{TCH} is a representation theorem for any
positive linear functional on $\mathbb{R}_{m}[t]\mid _{K}$ ($K\subseteq 
\mathbb{R}^{d}$ compact). We next formulate this result for the complex case
and give a new proof, based not on convexity but on $C^{\ast }$-algebra
ideas.

\begin{proposition}
\label{prop3.6}(cf. \cite[Th\'eor\`eme II]{TCH}) Let $m>0$ and let $K$ be a
compact subset of $\mathbb{C}^{d}$. If $\Phi :\mathbb{C}_{m,d}[z,\bar{z}%
]\mid _{K}\rightarrow \mathbb{C}$ is a positive linear functional, then
there exist $N\leq \dim \mathbb{C}_{m,d}[z,\bar{z}]\mid _{K}$, nodes $%
z_{1},\dots ,z_{N}$ in $K$, and positive weights $\rho _{1},\dots ,\rho _{N}$%
, such that $\Phi (p)=\sum_{i=1}^{n}\rho _{i}p(z_{i},\bar{z}_{i})$ $(p\in 
\mathbb{C}_{m,d}[z,\bar{z}]\mid _{K})$.
\end{proposition}

\begin{proof}
$\mathbb{C}_{m,d}[z,\bar{z}]$ is an \textit{operator system} in $C(K)$ (cf.\ 
\cite[Chapter 2]{PAU} and \cite[Chapter 5, Definition 33.1]{CON2}), so $\Phi 
$ can be extended to a positive linear functional $\tilde{\Phi}$ on $C(K)$ ( 
\cite[Exercise 2.10]{PAU} or \cite[Chapter 5, Proposition 33.2(c)]{CON2}. By
the Riesz Representation Theorem, there is a positive Borel measure $\mu $, $%
\limfunc{supp}\mu \subseteq K$, such that $\tilde{\Phi}(f)=\int_{K}f\,d\mu $ 
$(f\in C(K))$. The result now follows by applying Theorem \ref{theorem3.5}
to $\mu $.
\end{proof}

\section{Unique representing measures in TCMP\label{sec4unique}}

Concerning the analogue of Question \ref{question1.3} for the truncated
complex moment problem, consider the following possible properties of $%
\gamma \equiv \gamma ^{(2n)}$: \newline
(P1) $\gamma $ has a unique representing measure; \newline
(P2) $M(\gamma ):=\{\nu :\nu $ is a representing measure for $\gamma \}$ is
weak-$\ast $ closed in $B_{1}(C_{0}(\mathbb{C})^{\ast })$; \newline
(P3) $\gamma $ admits a finitely atomic representing measure.\newline
Results of the preceding sections show that (P1) $\Rightarrow $ (P2) $%
\Rightarrow $ (P3). In the sequel we establish (P1) in two basic cases of
TCMP.

Recall that $\gamma $ is \textit{flat} if $M(n)\equiv M(n)(\gamma )\geq 0$
and $\limfunc{rank}M(n)=\limfunc{rank}M(n-1)$. In \cite[Corollary 5.14]
{TCMP1} we proved that if $\gamma $ is flat, then $\gamma $ admits a unique
representing measure having moments of all orders (and this measure is $%
\limfunc{rank}M(n)$-atomic).

\begin{proposition}
\label{prop4.1}If $\gamma^{(2n)}$ is flat, then there exists a unique
representing measure, which is $\limfunc{rank}M(n)$-atomic.
\end{proposition}

\begin{proof}
Let $\mu $ be a representing measure for $\gamma ^{(2n)}$. Since 
\cite[Corollary 5.14]{TCMP1} implies that $\gamma $ has a unique
representing measure having moments of all orders, it suffices to establish
that $\mu $ has moments of all orders. We first consider moments of order $%
2n+1$. Suppose $i,j\geq 0$ and $i+j=n$; since $M(n)\geq 0$ and $\limfunc{rank%
}M(n)=\limfunc{rank}M(n-1)$, there exists $p_{ij}\in \mathbb{C}_{n-1}[z,\bar{%
z}]$ such that $\bar{Z}^{i}Z^{j}=p_{ij}(Z,\bar{Z})\in \mathcal{C}_{M(n)}$
(the column space of $M(n)$). Since $\mu $ is a representing measure, $\bar{z%
}^{i}z^{j}=p_{ij}(z,\bar{z})$ on $\limfunc{supp}\mu $ \cite[Proposition 3.1]
{TCMP1}. Thus, for $k,l\geq 0$, $k+l=n+1$, $\bar{z}^{i+k}z^{j+l}=(\bar{z}%
^{k}z^{l}p_{ij})$ $(z,\bar{z})$ on $\limfunc{supp}\mu $, and since $\deg 
\bar{z}^{k}z^{l}p_{ij}\leq 2n$, then $\int \bar{z}^{i+k}z^{j+l}\,d\mu $ is
convergent. By considering all indices $i$, $j$, $k$, $l$ with $i+j=n$ and $%
k+l=n+1$, it follows that $\mu $ has convergent moments up to degree $2n+1$.

We next consider degree $2n+2$. There exists $p_{0,n}\in\mathbb{C}_{n-1}[z,%
\bar{z}]$ such that $Z^{n}=p_{0,n}(Z,\bar{Z})$, whence $z^{n}=p_{0,n}(z,\bar{%
z})$ on $\limfunc{supp}\mu$. Thus $\left| z\right| ^{2n+2}=(\left| z\right|
^{2}\left| p_{0,n}\right| ^{2})(z,\bar{z})$ on $\limfunc{supp}\mu$, whence $%
\deg\left| zp_{0,n}\right| ^{2}\leq2n$ and $\int\left| z\right|
^{2n+2}\,d\mu<+\,\infty$. Thus, for all $(i,j)$ such that $i+j=2n+2$, $\int%
\bar{z}^{i}z^{j}\,d\mu$ is absolutely convergent, hence convergent.

Since $\mu $ has convergent moments up to degree $2n+2$, we may consider $%
M(n+1)[\mu ]$. Since $\mu $ is a representing measure for $M(n+1)[\mu ]$, it
follows that $M(n+1)[\mu ]$ is positive and recursively generated. For $%
i+j=n $, we have $\bar{Z}^{i}Z^{j}=p_{ij}(Z,\bar{Z})$ in $\mathcal{C}_{M(n)}$%
, and since $M(n+1)[\mu ]\geq 0$, it follows that $\bar{Z}^{i}Z^{j}=p_{ij}(Z,%
\bar{Z})$ in $\mathcal{C}_{M(n+1)[\mu ]}$ \cite[Proposition 2.4 (Extension
Principle)]{FIA1}. By recursiveness, we have $\bar{Z}^{i}Z^{j+1}=(zp_{ij})(Z,%
\bar{Z})$ and $\bar{Z}^{i+1}Z^{j}=(\bar{z}p_{ij})(Z,\bar{Z})$ in $\mathcal{C}%
_{M(n+1)[\mu ]}$. Since $\deg zp_{ij}\leq n$ and $\deg \bar{z}p_{ij}\leq n$,
it follows that for $k+l=n+1$, $\bar{Z}^{k}Z^{l}\in \left\langle \bar{Z}%
^{i}Z^{j}\right\rangle _{0\leq i+j\leq n}$ in $\mathcal{C}_{M(n+1)[\mu ]}$,
whence $\limfunc{rank}M(n+1)[\mu ]=\limfunc{rank}M(n)[\mu ]$.

Thus, $\gamma^{2(n+1)}[\mu]$ is flat. The preceding argument may be repeated
to produce successive flat extensions $M(n+1)[\mu]$, $M(n+2)[\mu],\dots$;
thus $\mu$ has convergent moments of all orders.
\end{proof}

\begin{proposition}
\label{prop4.2}If $\gamma ^{(2n)}$ has a representing measure and if $%
M(n)(\gamma )$ admits an analytic relation as in (\ref{equ13}) for some $%
k\leq n$, then $\gamma $ has a unique representing measure, which is
finitely atomic with at most $k^{2}$ atoms.
\end{proposition}

\begin{remark}
\cite[Theorem 3.1]{FIA2} illustrates the case of $\gamma \equiv \gamma
^{(6)} $ such that $M(3)(\gamma )$ has an analytic relation with $k=n=3$,
and such that the unique representing measure is $9$-atomic.
\end{remark}

We require the following preliminary result, which is of independent
interest.

\begin{proposition}
\label{lemma4.3}A polynomial of the form $p(z,\bar{z})\equiv z^{k}-q(z,\bar{z%
})$, where $q\in \mathbb{C}_{k-1}[z,\bar{z}]$, has at most $k^{2}$ roots.
\end{proposition}

\begin{proof}
The result is obvious for $k=1$, so we assume $k\geq 2$. Let $\mathcal{Z}(p)$
denote the zero set of $p$ and assume that $\mathcal{Z}(p)\supseteq \Lambda
\equiv \{z_{1},\dots ,z_{m}\}$, where $m:=k^{2}+1$ and the $z_{i}$s are
distinct. Let $\mu $ be a positive measure with $\limfunc{supp}\mu =\Lambda $%
, let $n:=2k-2$, and form $M(n)[\mu ]$. Clearly $\mu $ is a representing
measure for $M(n)[\mu ]$, and since $\deg p\leq n$ and $\limfunc{supp}\mu
\subseteq \mathcal{Z}(p)$, we must have $p(Z,\bar{Z})=0$ in $\mathcal{C}%
_{M(n)[\mu ]}$ (\cite[Proposition 3.1]{TCMP1}). It follows that $Z^{k}=q(Z,%
\bar{Z})$, with $k=\frac{n}{2}+1=\left[ \frac{n}{2}\right] +1$ and $\deg q<k$%
. By \cite[Theorem 3.1]{TCMP2}, $M(n)[\mu ]$ admits a unique positive,
recursively generated extension $M\equiv M(\infty )$, and $M$ is a flat
extension. Since $M(\infty )[\mu ]$ is positive and recursively generated,
we have $M(\infty )[\mu ]=M$, and \cite[Proposition 4.6]{TCMP1} implies that 
\begin{equation*}
m=\limfunc{card}\limfunc{supp}\mu =\limfunc{rank}M(\infty )[\mu ]=\limfunc{%
rank}M=\limfunc{rank}M(n)[\mu ].
\end{equation*}
Observe now that the size of $M(n)[\mu ]$ is $\frac{(n+1)(n+2)}{2}=k(2k-1)$.
The column relation $Z^{k}=q(Z,\bar{Z})$ (and the conjugate relation $\bar{Z}%
^{k}=\bar{q}(Z,\bar{Z})$) propagate recursively (and disjointly) in $%
\mathcal{C}_{M(n)[\mu ]}$ to generate $2(1+2+\cdots +k-1)=k(k-1)$ columns
each of which is a linear combination of columns of strictly lower degree.
(These $k(k-1)$ columns are $Z^{k}$, $\bar{Z}^{k}$, $Z^{k+1}$, $\bar{Z}Z^{k}$%
, ..., $\bar{Z}^{k}Z$, $\bar{Z}^{k+1}$,..., $Z^{n}$, ..., $\bar{Z}%
^{n-k}Z^{k} $, $\bar{Z}^{k}Z^{n-k}$, ..., $\bar{Z}^{n}$.) Thus, 
\begin{equation*}
m=\limfunc{rank}M(n)[\mu ]\leq k(2k-1)-k(k-1)=k^{2}<m,
\end{equation*}
a contradiction. We conclude that $\limfunc{card}\mathcal{Z}(p)\leq k^{2}$.
\end{proof}

\begin{remark}
\label{remark4.4}The following examples show that the estimate in
Proposition \ref{lemma4.3} is sharp. Let $p(z,\bar{z}):=z^{2}-\bar{z}$. Then 
$\mathcal{Z}(p)=\{0\}\cup \{\omega ^{j}:0\leq j\leq 2\}$, where $\omega $ is
a primitive cubic root of unity. Thus, $\limfunc{card}\mathcal{Z}(p)=4=(\deg
p)^{2}$. For $k=3,4,5$, A. Wilmshurst \cite[Chapter 3, Section 5]{WIL} has
constructed polynomials $q_{3}(z,\bar{z}):=z^{3}+z^{2}+z-2\bar{z}^{2}-2\bar{z%
}$, $q_{4}(z,\bar{z}):=z^{4}+3z^{2}-z-3\bar{z}^{3}-3\bar{z}^{2}$ and $%
q_{5}(z,\bar{z}):=z^{5}+5z^{3}-10z^{2}+5z+5\bar{z}^{4}-5\bar{z}^{3}$, with $%
\limfunc{card}\mathcal{Z}(q_{k})=(\deg q_{k})^{2}\;(k=3,4,5)$.
\end{remark}

\begin{proof}[Proof of Proposition $\mathrm{4.2}$]
We have an analytic relation of the form $Z^{k}=q(Z,\bar{Z})$ for some $%
k\leq n$, where $q\in \mathbb{C}_{k-1}[z,\bar{z}]$. Let $p(z,\bar{z}%
):=z^{k}-q(z,\bar{z})$; Proposition \ref{lemma4.3} implies that $\limfunc{%
card}\mathcal{Z}(p)\leq k^{2}$. Let $\mu $ be a representing measure for $%
\gamma $. Since $\limfunc{supp}\mu \subseteq \mathcal{Z}(p)$ 
\cite[Proposition 3.1]{TCMP1}, $\limfunc{card}\limfunc{supp}\mu \leq k^{2}$.
Thus $\mu $ has moments of all orders, and recursiveness implies that $%
M\equiv M(\infty )[\mu ]$ is completely determined by the analytic column
relations 
\begin{equation}
Z^{k+j}=(z^{j}q)(Z,\bar{Z})\qquad (j\geq 0).  \label{eq4.1}
\end{equation}
(Note that for each $p$, $M(p+1)[\mu ]$ is completely determined by $%
M(p)[\mu ]$ and the moments in $Z^{p+1}$.) Further, since $\mu $ is a
representing measure for $M$, $\limfunc{rank}M\leq \limfunc{card}\limfunc{%
supp}\mu \leq k^{2}$. Now, if $\nu $ is a representing measure for $\gamma $%
, then it follows as above that $M(\infty )[\nu ]$ is also determined by (%
\ref{eq4.1}), whence $M(\infty )[\nu ]=M$. Since $M$ is positive and has
finite rank, \cite[Theorem 4.7]{TCMP1} implies that it has a unique
representing measure, whence $\mu =\nu $.
\end{proof}

\end{document}